\documentclass[ a4paper]{amsart}%
\usepackage{amsmath,amsthm,amssymb,amsfonts}
\usepackage{enumerate,amscd,amsxtra}
\usepackage{mathrsfs}
\usepackage{graphicx}
\usepackage[cmtip,all]{xy}
\usepackage{amsmath}
\usepackage{amsfonts}
\usepackage{amssymb}%
\providecommand{\U}[1]{\protect\rule{.1in}{.1in}}
\nonstopmode
\newcounter{fig}

\theoremstyle{plain}
\newtheorem{Theorem}{Theorem}[section]

\newtheorem{corollary}[Theorem]{Corollary}

\newtheorem{lemma}[Theorem]{Lemma}

\newtheorem{proposition}[Theorem]{Proposition}
\newtheorem{theorem}[Theorem]{Theorem}

\newtheorem{thm}[Theorem]{Theorem}
\newtheorem{con}[Theorem]{Conjecture}
\theoremstyle{definition}

\newtheorem{remark}[Theorem]{Remark}

\theoremstyle{remark}
\def\Changed/{\ifvmode\else\vadjust{\vbox to 0pt{\vskip -\baselineskip\hbox to 0pt{\hss\vrule height 0pt depth 1.2\baselineskip\hskip 1em}\vss}}\fi}
\def\Math#1{\def\MathString{#1}\futurelet\MathDelim\MathChoose}
\def\MathChoose{\ifmmode\let\MathDo\MathString              \else\let\MathDo\MathSkip\fi              \MathDo}
\def\MathSkip{\ifx\MathDelim/\def\MathDo{$\MathString$\EatOne}              \else\def\MathDo{$\MathString$}\fi              \MathDo}
\def\Text#1{\def\TextString{#1}\futurelet\TextDelim\TextSkip}
\def\TextSkip{\ifx\TextDelim/\def\TextDo{\TextString\EatOne}              \else\let\TextDo\TextString\fi              \TextDo}
\def\EatOne#1{}
\def\SkipToEndScan#1\EndScan{}
\def\Scan#1#2#3{\ifx#1#2#3\expandafter\SkipToEndScan\fi\Scan#1}
\def\Upper#1{\Scan#1aAbBcCdDeEfFgGhHiIjJkKlLmMnNoOpPqQrRsStTuUvVwWxXyYzZ#1#1\EndScan}
\def\Phrase#1 #2/#3/#4=#5 #6/#7/#8.{\expandafter\edef\csname#2#3\endcsname{\noexpand\Text{#6#7}}
\expandafter\edef\csname\Upper#2#3\endcsname{\noexpand\Text{\Upper#6#7}}
\expandafter\edef\csname#1#2#3\endcsname{\noexpand\Text{#5 #6#7}}
\expandafter\edef\csname\Upper#1#2#3\endcsname{\noexpand\Text{\Upper#5 #6#7}}
\expandafter\edef\csname#2#4\endcsname{\noexpand\Text{#6#8}}
\expandafter\edef\csname\Upper#2#4\endcsname{\noexpand\Text{\Upper#6#8}}
}

\begin{document}
\title[Acyclic groups]{The acyclic group dichotomy}
\author{A. J. Berrick}
\address{Department of Mathematics, National University of Singapore, Kent Ridge
119076, SINGAPORE}
\email{berrick@math.nus.edu.sg}
\dedicatory{To the memory of Karl Gruenberg -- fine mathematician, exemplary colleague,
and dear friend.}
\begin{abstract}
Two extremal classes of acyclic groups are discussed. For an arbitrary group
$G$, there is always a homomorphism from an acyclic group of cohomological
dimension $2$ onto the maximum perfect subgroup of $G$, and there is always an
embedding of $G$ in a binate (hence acyclic) group. In the other direction,
there are no nontrivial homomorphisms from binate groups to groups of finite
cohomological dimension. Binate groups are shown to be of significance in
relation to a number of important $K$-theoretic isomorphism conjectures.

\end{abstract}
\maketitle

\section{Introduction}

The theme of this note is that the world of acyclic groups appears to be
dominated by two very distinct classes.

Recall that a discrete group $G$ is called \emph{acyclic}\textit{ }if all
ordinary homology groups, with trivial integer coefficients, are zero in
positive dimensions. In particular, considering dimension $1$, all such groups
are perfect. Among acyclic groups, we focus on classes that we may think of as
the small and large acyclic groups, in the sense of cohomological dimension.
In what follows, we prove that every group is sandwiched between these
classes, in that it receives maps from small acyclic groups and embeds in a
large acyclic, \emph{binate }group. Schematically:%
\begin{equation}%
\begin{array}
[c]{ccccc}%
\text{acyclic group of cd }2 &  &  &  & \\
&  &  &  & \\
{}^{\exists}\downarrow^{\mathrm{onto}} &  &  &  & \\
\text{perfect radical} & \hookrightarrow & \text{arbitrary group} &
\overset{\exists}{\hookrightarrow} & \text{binate group}%
\end{array}
\label{small to large}%
\end{equation}
We prove a rigidity result to the effect that the flow is irreversible,
inasmuch as there are no nontrivial maps from large acyclic groups to small
ones. Arising from this setup, large acyclic groups play a crucial role in
certain long-standing conjectures. The note concludes with brief lists of
examples of each class, for which the reader is referred to
\cite{Berrick_Steer chapter} for details.

\section{The small: acyclic groups of finite cohomological dimension}

Arguably, the first (nontrivial) acyclic groups in the literature, those of
Baumslag and Gruenberg \cite{BaumslagGruenberg}, are of this kind. They
consider a group $G$ generated by two elements subject to a single relation,
and show that if its abelianization $G_{\mathrm{ab}}$ is infinite cyclic and
the commutator subgroup $G^{\prime}$ is perfect, then $G^{\prime}$ is acyclic.
Since two-generator one-relator groups have cohomological dimension $2$, so
does the acyclic group $G^{\prime}$. Remarkably, a geometrically-defined
example was provided almost simultaneously: in \cite{Epstein}, Epstein
performs a variant of the grope construction \cite{Stanko}, \cite{Teichner},
whereby a doubly-punctured torus is added at each stage. In \cite{BerrickWong}%
, Epstein's group is shown to be acyclic of Baumslag-Gruenberg type: it is the
commutator subgroup of the group
\[
G=\left\langle x,y\mid x=\left[  x,yx^{-1}y^{-1}\right]  \left[
x,y^{-1}xy\right]  \right\rangle .
\]

In general, one can describe a grope \cite{DranRepovs} as the direct limit $L$
of a nested sequence of compact $2$-dimensional polyhedra
\[
L_{0}\overset{}{\hookrightarrow}L_{1}\hookrightarrow L_{2}\hookrightarrow
\cdots
\]
obtained as follows. Take $L_{0}$ as some $S_{g}$, an oriented compact surface
of positive genus $g$ from which an open disk has been deleted. To form
$L_{n+1}$ from $L_{n}$, for each loop $a$ in $L_{n}$ that generates the group
$H_{1}(L_{n})$, attach to $L_{n}$ some $S_{g_{a}}$ by identifying the boundary
of $S_{g_{a}}$ with the loop $a$. Since the fundamental group of $S_{g}$
punctured is a free group on $2g$ generators, this procedure embeds each
$\pi_{1}(L_{n})$, and thus each finitely generated subgroup of $\pi_{1}(L)$,
as a subgroup of a free group, with each generator $a$ of $\pi_{1}(L_{n})$
becoming a product of $g_{a}$ commutators in $\pi_{1}(L_{n+1})$. Hence
$\pi_{1}(L)$ is a countable, perfect, locally free group. Since free groups
have geometric dimension $1$, and homology groups commute with direct limits,
there is a key (well-known) observation that we record.

\begin{lemma}
\label{perfect locally free is acyclic}Perfect, locally free groups are
acyclic, of geometric dimension $2$.$\hfill\Box\smallskip$
\end{lemma}

The converse to this lemma fails: \cite[p.13]{BDH} exhibits an acyclic group
of geometric dimension $2$ that is not locally free, indeed it contains the
fundamental groups of all closed orientable surfaces of genus $\,\geq2$.
Moreover, there are acyclic groups of all geometric dimensions, as in the next
result. For a stronger statement, we can speak of the \emph{finite geometric
dimension }of a group $G$, namely the minimum dimension of a finite CW-complex
that is a $K(G,1)$ (where such a complex exists).

\begin{theorem}
\cite{BerrickHillman} For any integer $n>1$, there is an acyclic group of
finite geometric dimension $n$.
\end{theorem}

The most economical example of $L$ as above is the \emph{minimal
grope}\textit{ }$M^{\ast}$, for which one takes only one genus one surface at
each step. Homotopically, each $L_{n}$ is in this case just a bouquet of
finitely many circles. So $M^{\ast}$ is the classifying space of the group
having: as generators, symbols $x_{w}$ for each $w\in\Sigma$, the set of all
(nonempty) words of finite length on the two symbols $0,1$; and as relations,
$x_{w}=[x_{w0},x_{w1}]$. Thus, the $x_{w}$ with $w$ of length $n$ generate the
free group $\pi_{1}(L_{n})$, but each is a commutator when embedded in the
larger free group $\pi_{1}(L_{n+1})$. For the $3$-dimensional version of this
construction, see \cite[Theorem 2]{FreedmanFreedman}, where a nested sequence
of handlebodies leads to a locally free fundamental group, and conversely.

To see the process more algebraically \cite{BerrCasa}, one takes the free
product of countable, perfect, locally free groups $F_{\mathbf{n}}$, one for
each sequence\textrm{ }$\mathbf{n}=(n_{1},n_{2},n_{3},\ldots)$ of positive
integers. $F_{\mathbf{n}}$ is the direct limit of the following direct system
$(F_{\mathbf{n},r},\varphi_{r})$. For $r=0$, define $F_{\mathbf{n},0}$ to be
infinite cyclic with a generator~$x_{0}$, and for each $r\geq1$, let
$F_{\mathbf{n},r}$ be the free group freely generated by the set of
$2^{r}n_{1}\cdots n_{r}$ symbols
\[
\left\{  x_{r}(\varepsilon_{1},\ldots,\varepsilon_{r};\,i_{1},\ldots
,i_{r})\mid\varepsilon_{k}\in\{0,1\},\;1\leq i_{k}\leq n_{k}\right\}  .
\]
Then for $r\geq0$, $\varphi_{r}\hspace{-0.1cm}:F_{\mathbf{n},r}\rightarrow
F_{\mathbf{n},r+1}$ maps $x_{r}(\varepsilon_{1},\dots,\varepsilon_{r}%
;\,i_{1},\ldots,i_{r})$ to the product of commutators%
\[
\prod_{i_{r+1}=1}^{n_{r+1}}[x_{r+1}(\varepsilon_{1},\ldots,\varepsilon
_{r},0;\,i_{1},\ldots,i_{r},i_{r+1}),\;x_{r+1}(\varepsilon_{1},\ldots
,\varepsilon_{r},1;\,i_{1},\ldots,i_{r},i_{r+1})].
\]

Now, for every perfect group $P$ and each element $x\in P$, one can write $x$
as a product of commutators. In consequence there is a sequence of integers
$n$ (which may be chosen to be increasing) and a homomorphism $\psi
:F_{\mathbf{n}}\rightarrow P$ whose image contains~$x$ \cite{Heller}. A free
product of such $F_{\mathbf{n}}$ then gives a perfect, locally free group with
image containing any finitely generated subgroup of $P$. Since $P$ is the
direct limit of its finitely generated subgroups, by taking the direct limit
of free products of groups $F_{\mathbf{n}}$ we obtain a perfect, locally free
group that has $P$ as a homomorphic image. Starting with an arbitrary group,
recall that the group generated by all its perfect subgroups is also perfect,
its \emph{perfect radical}. By the above procedure, this perfect radical is
the homomorphic image of a perfect, locally free group. The conclusion is thus
the left hand side of the schematic flow of diagram (\ref{small to large}), as follows.

\begin{proposition}
For any group $G$, there is a homomorphism $\phi:A\rightarrow G$, where $A$ is
a locally free acyclic group, such that every perfect subgroup of $G$ lies in
the image of $\phi$.\hfill$\Box$\smallskip
\end{proposition}

\section{The large: binate groups}

We begin with an elementary observation, that follows readily from the formula%
\[
\left[  u,\,ab\right]  =\left[  u,a\right]  \left[  u,b\right]  \left[
\left[  b,u\right]  ,\,a\right]  \text{,}%
\]
where $[u,v]=uvu^{-1}v^{-1}$.

\begin{lemma}
Let $H$ be a subgroup of a group $G$, $u\in G$, and $\varphi:H\rightarrow G$ a
set function such that for all $h\in H$,
\[
h=\left[  u,\,\varphi(h)\right]  .\label{binate equation copy(1)}%
\]
Then $\varphi$ is a homomorphism if and only if both $\varphi$ is injective
and $\left[  H,\,\varphi(H)\right]  =1$.\hfill$\Box$\smallskip
\end{lemma}

A nontrivial group $G$ is called \emph{binate}\textit{ }\cite{Berr sgtc} if
for any finitely generated subgroup $H$ of $G$ there exists a homomorphism
(called a \emph{structure map}) $\varphi=\varphi_{H}:H\rightarrow G$ and
element (called a \emph{structure element}) $u=u_{H}\in G$ such that for all
$h\in H$,
\[
h=\left[  u,\,\varphi(h)\right]  .\label{binate equation}%
\]
Thus, the finitely generated subgroups of $G$ are indeed arranged in commuting
isomorphic pairs, as the name \textquotedblleft binate\textquotedblright\ is
intended to suggest.

Evidently, a binate group is perfect; indeed, every element is a commutator.
Observe too that $u\in H$ would imply that $H=1$. Therefore, binate groups
cannot be finitely generated. Deeper facts are summarized as follows.

\begin{theorem}
\label{binate acyclic no repns}\emph{(a)} \textrm{\cite{Berr sgtc}} Binate
groups are acyclic.

\emph{(b)} \textrm{\cite{AlperinBerrick} }Binate groups have no
finite-dimensional representations over any field.
\end{theorem}

M. V. Sapir raised with me the interesting question of whether the pairing of
$H$ with $\varphi(H)$ above yields a copy of $H\times H$ in the binate group
$G$. While a weakened version of this query is readily answered (in (a)
below), in its original form the question is more subtle, because the centre
$\mathcal{Z}(H)$ of a finitely generated group $H$ need not be finitely generated.

\begin{theorem}
\label{01203a}\emph{(a)}\textrm{ }If $H$ is a finitely generated subgroup of a
binate group $G$, then there are structure maps $\varphi_{1},\varphi
_{2},\ldots$ for $H$ such that, for any distinct $i,j$, we have $\varphi
_{i}(H)\times\varphi_{j}(H)\leq G$. Thus, $G$ contains an infinite product of
copies of $H$.

\emph{(b)}\textrm{ }Let $H\leq H_{1}$ be finitely generated subgroups of a
binate group $G$ such that $H\cap\mathcal{Z}(H_{1})$ is finitely generated.
Then there exists a structure map $\psi=\psi_{H}$ for $H$ such that $H\cap
\psi(H)=1$; in other words, $H\times\psi(H)<G$.
\end{theorem}

\noindent\textbf{Proof. (a) }Write $H_{0}=H$ and, for $i\in\mathbb{N}$,
\[
H_{i}=\left\langle H_{i-1},\,\varphi_{i-1}(H_{i-1}),\,u_{i-1}\right\rangle
\]
where $\varphi_{i-1}$ is a structure map for $H_{i-1}$, with associated
structure element $u_{i-1}$. Thus $\left[  \varphi_{i}(H),\,\varphi
_{j}(H)\right]  =1$. If for some $i<j$ we have $a,b\in H$ with $\varphi
_{i}(a)=\varphi_{j}(b)$, then%
\[
a=\left[  u_{i},\,\varphi_{i}(a)\right]  =\left[  u_{i},\,\varphi
_{j}(b)\right]  =1,
\]
whence $\varphi_{i}(a)=1$. For the final assertion, recall that each
$\varphi_{i}$ is a monomorphism, and iterate.

\textbf{(b) }Since $H\cap\mathcal{Z}(H_{1})$ is a finitely generated abelian
group (hence isomorphic to $\mathbb{Z}^{\mathrm{rk}}\oplus\mathrm{Tor}$ with
$\mathrm{Tor}$ finite), we argue by induction on the pair
\[
(\mathrm{rk}_{\mathbb{Z}}(H\cap\mathcal{Z}(H_{1})),\,\left\vert \mathrm{Tor}%
(H\cap\mathcal{Z}(H_{1}))\right\vert ),
\]
ordered lexicographically. When this pair is $(0,1)$, we already have that
\[
H\cap\varphi_{H_{1}}(H)\leq H\cap\mathcal{Z}(H_{1})=1,
\]
so that the result is immediate with $\psi_{H}=\varphi_{H_{1}}$. For the
induction step, assuming the result for lesser values of the pair, suppose
that $x\in H-\{1\}$ has $\varphi(x)\in H$ where $\varphi=\varphi_{H_{1}}$.
Then $\varphi(x)\in H\cap\mathcal{Z}(H_{1})$, say $\varphi(x)=t^{n}$ for some
element $t$ of a minimal generating set for $H\cap\mathcal{Z}(H_{1})$. Since%
\[
1\neq x=\left[  u,t^{n}\right]  ,
\]
where $u=u_{H_{1}}$, and $\varphi$ is a monomorphism, one of the following
must occur.

\begin{enumerate}
\item[(i)] $t$ has infinite order. Then $x$ also must have infinite order.
Therefore $\left\langle t\right\rangle \cap\mathcal{Z}(\left\langle
H_{1},u\right\rangle )=1$. This implies that
\[
\mathrm{rk}_{\mathbb{Z}}(H\cap\mathcal{Z}(\left\langle H_{1},u\right\rangle
)<\mathrm{rk}_{\mathbb{Z}}(H\cap\mathcal{Z}(H_{1})).
\]

\item[(ii)] If
\[
\mathrm{rk}_{\mathbb{Z}}(H\cap\mathcal{Z}(\left\langle H_{1},u\right\rangle
)=\mathrm{rk}_{\mathbb{Z}}(H\cap\mathcal{Z}(H_{1})),
\]
then, since $x$ is nontrivial, $t$ must have finite order with $\left\langle
t\right\rangle \nsubseteq\mathcal{Z}(\left\langle H_{1},u\right\rangle )$.
Therefore $\,$%
\[
\left\vert \mathrm{Tor}(H\cap\mathcal{Z}(\left\langle H_{1},u\right\rangle
))\right\vert <\left\vert \mathrm{Tor}(H\cap\mathcal{Z}(H_{1}))\right\vert .
\]

\end{enumerate}

In either event, we can apply the induction hypothesis to obtain $\psi
_{H}=\psi_{\left\langle H_{1},u\right\rangle }$ as the desired structure map.
\hfill%
%TCIMACRO{\TeXButton{End Proof}{\quad\lower0.05cm\hbox{$\Box$}}}%
%BeginExpansion
\quad\lower0.05cm\hbox{$\Box$}%
%EndExpansion
\smallskip

\begin{remark}
This argument extends to include the case where $H\cap\mathcal{Z}(H_{1})$ is
the direct sum of a finitely generated group with a finite number of copies of
the rationals and a finite number of quasicyclic groups.
\end{remark}

An easy consequence of this result (to be strengthened below) is that our
classes of small and large acyclic groups are indeed disjoint.

\begin{corollary}
\label{01203b}Every binate group has infinite cohomological dimension.
\end{corollary}

\noindent\textbf{Proof. }By the theorem applied initially to a cyclic
subgroup, a binate group contains either free abelian groups of arbitrarily
high rank, or elementary abelian groups of arbitrarily high rank. However,
neither of these possibilities can occur in a group of finite cohomological
dimension. \hfill%
%TCIMACRO{\TeXButton{End Proof}{\quad\lower0.05cm\hbox{$\Box$}}}%
%BeginExpansion
\quad\lower0.05cm\hbox{$\Box$}%
%EndExpansion
\smallskip

Here is another result, due to I. Agol (private communication), that shows the
aversion of low-dimensional geometry for binate groups.

\begin{theorem}
\label{01131a copy(1)}$3$-manifold groups cannot be binate.
\end{theorem}

\noindent\textbf{Proof. }The important property is that for a binate group
$G$, and any finitely generated subgroup $H$, there is an isomorphic subgroup
$H^{\prime}$ such that $[H,H^{\prime}]=1$. First, we observe that $G$ cannot
be a free product, since if $G=G_{1}\ast G_{2}$, then taking a finitely
generated subgroup of the form $H=\left\langle h_{1},h_{2}\right\rangle $,
with $h_{1}$ in $G_{1}-\{1\}$, $h_{2}$ in $G_{2}-\{1\}$, it's clear that there
can be no subgroup $H^{\prime}<G$ commuting with $H$. This can be seen
geometrically, \textsl{e.g.} take a $K(G_{1},1)$ and $K(G_{2},1)$, and link
them together with an interval to get a $K(G,1)$. If we have a commutator
$[h_{1},h^{\prime}]=1$, then we get a map of a torus into the $K(G,1)$. This
torus map can be deformed disjoint from the interval, otherwise some loop in
the torus would be trivial, in which case some product $h_{1}^{a}$ $h^{\prime
b}$ would be trivial, which would imply that $h^{\prime}$ is in $G_{1}$, and
couldn't commute with $h_{2}$. Thus, if we deformed the torus into
$K(G_{1},1)$, we would get a conjugacy of $h^{\prime}$ into $G_{1}$.
Similarly, if $h^{\prime}$ commutes with $h_{2}$, then it is conjugate into
$H_{2}$, which would be a contradiction. So we may assume that $G$ is not a
free product. Thus, if $G=\pi_{1}(M^{3})$, then we may assume that $M^{3}$ is
irreducible and noncompact (since $G$ is not finitely generated), and
therefore a $K(G,1)$. However, by Corollary \ref{01203b}, $G$ cannot have
finite geometric dimension.

Alternatively, arguing geometrically, we can apply the Scott core theorem to
deduce that any finitely generated subgroup of $G$ is the fundamental group of
a compact irreducible $3$-manifold with nontrivial boundary. So, if we take
finitely generated non-virtually-cyclic subgroups $H$ and $H^{\prime}$ such
that $[H,H^{\prime}]=1$, then we get a finitely generated subgroup
$\left\langle H,H^{\prime}\right\rangle $ which is the fundamental group of a
compact manifold $N$ with boundary. Each pair of commuting elements in $H$ and
$H^{\prime}$ gives a $\pi_{1}$-injective map of a torus into $N$. By the
characteristic submanifold theorem, there is a canonical decomposition of $N$
along embedded incompressible tori, such that the complementary pieces are
simple or Seifert fibered, and such that any map of a torus into $N$ may be
homotoped into one of the Seifert fibered pieces. Non-closed Seifert fibered
spaces have a finite-sheeted covering which is of the form surface$\,\times
S^{1}$, so it has fundamental group of form free group$\,\times\mathbb{Z}$.
It's clear that the only pairs of commuting groups in here are of the form
$\mathbb{Z\times Z}$. Thus, $H$ must be a virtually cyclic group, a
contradiction. \hfill%
%TCIMACRO{\TeXButton{End Proof}{\quad\lower0.05cm\hbox{$\Box$}}}%
%BeginExpansion
\quad\lower0.05cm\hbox{$\Box$}%
%EndExpansion
\medskip

This result should be seen alongside the result of \cite{BerrickHillman} that
nontrivial acyclic $3$-manifold groups cannot be finitely generated. For
examples of such groups, see Appendix A below.

We turn now to the right-hand half of the flow diagram (\ref{small to large}).
A universal method for embedding any group $H$ in a binate group is the
\emph{universal binate tower }on $H$, constructed in \cite{Berr sgtc} by means
of HNN-extensions, as follows. Let $H_{0}=H$ and for each $i\geq0$
\[
H_{i+1}=\operatorname{gp}\left\langle H_{i}\times H_{i},u_{i}\mid
(g,g)=u_{i}(1,g)u_{i}^{-1}\hspace{0.15in}\text{for each }g\in H_{i}%
\right\rangle \text{,}%
\]
with $H_{i}$ embedded in $H_{i+1}$ as $H_{i}\times1$ and $\varphi
_{i}(g)=(1,g)$. The direct limit of the $H_{i}$ is a binate group (even though
the homology of the groups $H_{i}$ grows hyperexponentially with $i$). This
tower is the initial object in a category of binate towers with base group $H$
\cite{BerrVarad}.

\begin{theorem}
\textrm{\cite{BerrVarad} \emph{(a)} }Every group $H$ embeds in a universal
binate group.

\emph{(b)}\textrm{ }Every binate group that contains $H$ contains infinitely
many images of the universal binate group containing $H$.
\end{theorem}

Observe that, because homology cycles are always supported on finitely
generated subgroups, and because homology commutes with direct limits, in an
acyclic group $G$ every countable subgroup lies in a countable acyclic
subgroup of $G$. A similar result holds for binate groups $G$, by (b) above.

\begin{corollary}
Every countable subgroup of a binate group $G$ lies in a countable binate
subgroup of $G$.\hfill$\Box$\smallskip
\end{corollary}

An alternative way of embedding a given group in a binate group is provided by
the Kan-Thurston construction of the cone of a group \cite{KanThurston},
discussed in Appendix B below.

We are now in a position to prove the rigidity result stated in the Introduction.

\begin{theorem}
\label{01213a}Every homomorphism from a binate group to a group that is
residually of finite virtual cohomological dimension is trivial.
\end{theorem}

\noindent\textbf{Proof. }First, if every map from a group $G$ to a member of a
class of groups $\mathcal{X}$ is trivial, then every map from $G$ to a group
that is residually in $\mathcal{X}$ must also be trivial, because any
nontrivial element in the image of $G$ would be mapped nontrivially to a
member of $\mathcal{X}$. Next, by Theorem \ref{binate acyclic no repns}%
\thinspace(b) binate groups have no nontrivial finite quotients; thus, we may
ignore the word \textquotedblleft virtual\textquotedblright.

Since subgroups of groups of finite cohomological dimension also have finite
cohomological dimension and so are torsionfree, it suffices to show that any
torsionfree quotient $Q$ of a binate group $G$ contains a free abelian
subgroup of infinite rank. To do this, let $\pi:G\twoheadrightarrow Q$ be an
epimorphism and $x\in G-\mathrm{Ker}\pi$. Define $H_{0}=\left\langle
x\right\rangle $ and, for $i\in\mathbb{N}$,
\[
H_{i}=\left\langle H_{i-1},\,\varphi_{H_{i-1}}(H_{i-1}),\,u_{i-1}\right\rangle
\text{.}%
\]
Then the subgroup of $Q$ generated by all $\pi\varphi_{H_{i}}(H_{0})$ is
abelian. Its rank must be infinite because if for some $i<j$ and $m,n$ we have
$\pi\varphi_{H_{i}}(x)^{m}=\pi\varphi_{Hj}(x)^{n}$, then%
\[
\pi(x)^{m}=\left[  \pi(u_{i}),\,\pi\varphi_{H_{i}}(x^{m})\right]  =\pi(\left[
u_{i},\,\varphi_{H_{j}}(x^{n})\right]  )=1\text{,}%
\]
leaving $m=n=0$ as the only possibility. \hfill%
%TCIMACRO{\TeXButton{End Proof}{\quad\lower0.05cm\hbox{$\Box$}}}%
%BeginExpansion
\quad\lower0.05cm\hbox{$\Box$}%
%EndExpansion

\section{Binate groups as conjecture-testers}

The first hint that binate groups form an important class of groups for
testing some famous conjectures comes from \cite{BerrVarad}, where links to
the Kervaire conjecture are discussed. Recall that this conjecture asserts
that for any group $G$ either $G$ is trivial or the free product
$G\ast\mathbb{Z}$ with an infinite cyclic group cannot be normally generated
by a single element. By considering abelianizations one immediately sees that
any counterexample to the conjecture must be perfect.

\begin{Theorem}
\emph{\cite{BerrVarad} }The Kervaire conjecture holds for all groups if and
only if it holds for all binate groups.
\end{Theorem}

Remarkably, \cite{BerrVarad} also shows that in any binate group the Kervaire
conjecture is locally true, in that every finitely generated subgroup is
contained in one with the property of the conjecture.

Next, we consider the Bass trace conjecture \cite{Bass: InventM}, as follows.
With $\mathbb{Z}G$ the integral group ring of a group $G$, the
\emph{augmentation trace} is the $\mathbb{Z}$-linear map
\[
\epsilon:\mathbb{Z}G\rightarrow\mathbb{Z},\quad g\mapsto1\,
\]
induced by the trivial group homomorphism on $G$. The Hochschild homology
group $HH_{0}(\mathbb{Z}G)=\mathbb{Z}G/[\mathbb{Z}G,\mathbb{Z}G]$, with
$[\mathbb{Z}G,\mathbb{Z}G]$ the additive subgroup of $\mathbb{Z}G$ generated
by the elements $gh-hg$ ($g,h\in G$), identifies with $\bigoplus
_{[s]\in\lbrack G]}\mathbb{Z}\cdot\lbrack s]$, where $[G]$ is the set of
conjugacy classes $[s]$ of elements $s$ of $G$. The \emph{Hattori-Stallings
trace} of $M=\sum_{g\in G}m_{g}g\in\mathbb{Z}G$ is then defined by
\begin{align*}
\mathrm{HS}(M)  &  =M+[\mathbb{Z}G,\mathbb{Z}G]\\
&  =\sum_{[s]\in\lbrack G]}\epsilon_{s}(M)[s]\ \in\bigoplus_{\lbrack
s]\in\lbrack G]}\mathbb{Z}\cdot\lbrack s],
\end{align*}
where for $[s]\in\lbrack G]$, $\epsilon_{s}(M)=\sum_{g\in\lbrack s]}m_{g}$ is
a partial augmentation. Now, an element of $K_{0}(\mathbb{Z}G)$ is represented
by a difference of finitely generated projective $\mathbb{Z}G$-modules, each
of which is determined by an idempotent matrix having entries in $\mathbb{Z}%
G$. Combining the usual trace map to $\mathbb{Z}G$ of such a matrix with the
above trace on $\mathbb{Z}G$ induces the Hattori-Stallings trace homomorphism
on $K_{0}(\mathbb{Z}G)$, which is natural with respect to group homomorphisms.

\begin{con}
\label{Bass} For any group $G$, the induced map
\[
\mathrm{HS}_{G}:K_{0}(\mathbb{Z}G)\longrightarrow HH_{0}(\mathbb{Z}%
G)=\bigoplus_{[s]\in\lbrack G]}\mathbb{Z}\cdot\lbrack s]
\]
has image in $\mathbb{Z}\cdot\lbrack e]$.
\end{con}

There is likewise a version of the conjecture for the complex group ring, as
well as for intermediate rings. For a recent description of groups $G$
satisfying the conjecture, see \cite{BCM}, where the following adaptation of
the universal binate group is introduced. With the notation
\begin{align*}
\Delta_{C}  &  =\left\{  (c,c)\in C\times C\mid c\in C\right\}  ,\\
\Delta_{F}^{\prime}  &  =\left\{  (f^{-1},f)\in F\times F\mid f\in F\right\}
,
\end{align*}
let $G$ be a group with $\{F_{i}\}_{i\in I}$ as the set of all its finitely
generated abelian subgroups. For $i\in I$, write $C_{i}$ for the centralizer
in $G$ of $F_{i}$. Then in $G\times G$ the subgroups $(1\times F_{i}%
)=\{(1,f)\mid f\in F_{i}\}$ and $\Delta_{C_{i}}$ commute, so that their
product $(1\times F_{i})\cdot\Delta_{C_{i}}$ is also a subgroup of $G\times
G$. Likewise, $\Delta_{F_{i}}^{\prime}\cdot(1\times C_{i})$ is also a
subgroup, and the obvious bijection
\begin{align*}
(1\times F_{i})\cdot\Delta_{C_{i}}  &  \longrightarrow\Delta_{F_{i}}^{\prime
}\cdot(1\times C_{i})\\
(k,fk)  &  \longmapsto(f^{-1},fk)
\end{align*}
is a group isomorphism. Now define $A_{1}(G)$ to be the generalized HNN
extension
\[
A_{1}(G)=\mathrm{HNN}(G\times G;\ (1\times F_{i})\cdot\Delta_{C_{i}}%
\cong\Delta_{F_{i}}^{\prime}\cdot(1\times C_{i}),\ t_{i})_{i\in I}%
\]
meaning that, whenever $f\in F_{i}$ and $k\in C_{i}$,
\[
(k,fk)=t_{i}(f^{-1},fk)t_{i}^{-1}.
\]
For $n\geq2$, inductively define $A_{n}(G)=A_{1}(A_{n-1}(G))$. Since
$A_{n-1}(G)\leq A_{n}(G)$, we put $A(G)=\mathbb{\cup}A_{n}(G)$.

\begin{thm}
\label{Berrick}\cite{BCM} The homomorphism $G\rightarrow A(G)$ has the
following properties.

\begin{enumerate}
\item[(a)] It is a functorial Frattini embedding.

\item[(b)] Every finitely generated abelian subgroup of $A(G)$ has its
centralizer in $A(G)$ binate.

\item[(c)] The prime powers that occur as orders of elements of $A(G)$ are
precisely those that occur as orders of elements of $G$.

\item[(d)] If $G$ is infinite, then $A(G)$ has the same cardinality as $G$.
\end{enumerate}
\end{thm}

Here, recall that a Frattini embedding is one in which nonconjugate elements
remain nonconjugate. (Also, observe that by taking the trivial subgroup in
(b), one sees that $A(G)$ itself is binate.) In \cite{BCM}, properties (a) and
(b) are combined, with (b) used to show that the equivariant $K$-homology
$K_{\ast}^{A}($\underline{$E$}$A)$ of the universal proper $A$-CW-complex
\underline{$E$}$A$ reduces to a direct sum of complex vector spaces afforded
by the complex representation rings of the finite subgroups of $A$. Since the
Hattori-Stallings trace $\mathrm{HS}_{A}$ is known to exclude nontrivial
finite conjugacy classes from its image \cite{BerrHesselholt}, \cite{Linnell},
it follows that it yields only the trivial conjugacy class on elements of
$K_{0}(\mathbb{C}A)$ that originate in $K_{0}^{A}($\underline{$E$}$A)$. When
$G$ satisfies the Bost conjecture \cite{Lafforgue}, then it can be deduced
that the image of the Hattori-Stallings trace $\mathrm{HS}_{G}$ in
$HH_{0}(\mathbb{C}A)$ is just the trivial conjugacy class. Thereby, using (a),
\cite{BCM} shows that if $G$ satisfies the Bost conjecture, then it also
satisfies the Bass conjecture.

The point of (a) of the theorem is that for any ring $R$ the induced map
$HH_{0}(RG)\rightarrow HH_{0}(RA)$ is a monomorphism. Therefore, from the
commuting square%
\[%
\begin{array}
[c]{ccc}%
K_{0}(\mathbb{Z}G) & \longrightarrow & K_{0}(\mathbb{Z}A)\\
\downarrow^{\mathrm{HS}_{G}} &  & \downarrow^{\mathrm{HS}_{A}}\\
HH_{0}(\mathbb{Z}G) & \rightarrowtail & HH_{0}(\mathbb{Z}A)
\end{array}
\]
we immediately obtain a noteworthy consequence.

\begin{corollary}
The Bass conjecture holds for all groups if and only if it holds for all
binate groups.$\hfill\Box\smallskip$
\end{corollary}

As noted in \cite{Mislin}, there are several related conjectures (mostly for
torsionfree groups), such as the weak Bass conjecture, Atiyah conjecture,
trace conjecture, strong trace conjecture, embedding conjecture and
zerodivisor conjecture, as well as the fibered Farrell-Jones conjecture
\cite{BartelsLueck}, that are subgroup-closed. It is therefore immediate from
Theorem \ref{Berrick}\thinspace(a) that they hold for all (respectively all
torsionfree) groups if and only if they hold for all (resp. all torsionfree)
binate groups. With regard to the weak Bass conjecture, there is a further
simplification available. For, according to Bass (see \cite{BCM Htpy
idempotents}), this conjecture holds for all groups if and only if it is valid
for all finitely presented groups. Now, from \cite[Theorem 6]{BeMi}, there is
a finitely presented, strongly torsion generated, acyclic group, $\mathrm{BM}$
say, with the property that it contains (an isomorphic copy of) every finitely
presented group. These two simplifications thereby combine to reveal that a
single binate group suffices to test the validity of the conjecture.

\begin{corollary}
The following assertions are equivalent:

\begin{enumerate}
\item[(i)] the weak Bass conjecture holds for all groups;

\item[(ii)] the weak Bass conjecture holds for the finitely presented acyclic
group $\mathrm{BM}$;

\item[(iii)] the weak Bass conjecture holds for the universal binate group on
$\mathrm{BM}$;

\item[(iv)] the weak Bass conjecture holds for the binate group $A(\mathrm{BM}%
)$.$\hfill\Box\smallskip$
\end{enumerate}
\end{corollary}

Property (c) of the theorem sheds light on two further conjectures. With the
notation that $i$ is the map induced on $K$-theory from the standard inclusion
of the reduced $C^{\ast}$-algebra in the group von Neumann algebra, and
$\mathrm{tr}_{\mathcal{N}(G)}$ is the map induced by the standard
complex-valued von Neumann trace on $\mathcal{N}(G)$, the modified trace
conjecture of Baum-Connes \cite{BaumConnes}, \cite[Conjecture 0.2]{Lueck
Baum-Connes trace conjecture} asserts that, for any discrete group $G$, the
image of the Kaplansky trace%
\[
\kappa:K_{0}(C_{r}^{\ast}(G))\overset{i}{\longrightarrow}K_{0}(\mathcal{N}%
(G))\overset{\mathrm{tr}_{\mathcal{N}(G)}}{\longrightarrow}\mathbb{C}%
\]
lies in the ring $\Lambda_{G}$ of rational numbers whose denominators are
products of the orders of finite subgroups of $G$. We note that, as in
\cite[p.616]{BCM}, \cite[p.\thinspace58]{Mislin}, orders of finite subgroups
do indeed occur as denominators.

Also, there is the Farkas conjecture \cite{Farkas_Group rings questionnaire
Comm Alg 80}, which asserts that a discrete group $G$ has elements of all
prime orders that occur as denominators of nonintegral rational numbers in the
image of the composition%
\[
K_{0}(\mathbb{C}G)\longrightarrow K_{0}(C_{r}^{\ast}(G))\overset
{i}{\longrightarrow}K_{0}(\mathcal{N}(G))\overset{\mathrm{tr}_{\mathcal{N}%
(G)}}{\longrightarrow}\mathbb{C}\text{.}%
\]
Now, observe from (c) of Theorem \ref{Berrick} that $\Lambda_{G}%
=\Lambda_{A(G)}$, and that if $A(G)$ has an element of order $p$, then so does
$G$. On the other hand, since $G$ embeds in $A(G)$, there is a commuting
diagram%
\[%
\begin{array}
[c]{ccccccc}%
K_{0}(\mathbb{C}G) & \longrightarrow & K_{0}(C_{r}^{\ast}(G)) & \overset
{i}{\longrightarrow} & K_{0}(\mathcal{N}(G)) &  & \\
&  &  &  &  & \overset{\mathrm{tr}_{\mathcal{N}(G)}}{\searrow} & \\
\downarrow &  & \downarrow &  & \downarrow &  & \mathbb{C}\\
&  &  &  &  & \overset{\mathrm{tr}_{\mathcal{N}(A(G))}}{\nearrow} & \\
K_{0}(\mathbb{C}A(G)) & \longrightarrow & K_{0}(C_{r}^{\ast}(A(G))) &
\overset{i}{\longrightarrow} & K_{0}(\mathcal{N}(A(G))) &  &
\end{array}
\]
By combining these facts, we immediately obtain the following.

\begin{corollary}
\emph{(a)}\textrm{ }The modified trace conjecture of Baum-Connes holds for all
groups if and only if it holds for all binate groups.

\emph{(b)}\textrm{ }The Farkas conjecture holds for all groups if and only if
it holds for all binate groups.$\hfill\Box\smallskip$
\end{corollary}

In this connection, we note that for binate groups $A$ there is a stronger
possible result. For, if $F_{\alpha}$ and $F_{\beta}$ are finite subgroups of
a binate group $A$, then together they generate a finitely generated subgroup
$H$ of $A$. Thus, the structure map $\varphi_{H}$ yields by Theorem
\ref{01203a} above a copy of $H\times H$, which contains $F_{\alpha}\times
F_{\beta}$ as a subgroup. Thereby, $A$ contains finite subgroups isomorphic to
any direct product of the form $F_{\alpha}\times\cdots\times F_{\alpha}\times
F_{\beta}\times\cdots\times F_{\beta}$. Hence, the whole ring $\Lambda_{A}$
occurs in the image of $\kappa$.

\begin{proposition}
If the modified trace conjecture of Baum-Connes holds for a binate group $A$,
then the image of its Kaplansky trace is precisely $\Lambda_{A}$.$\hfill
\Box\smallskip$
\end{proposition}

Finally, we remark that the proof of Atiyah's $L^{2}$-index theorem in
\cite{ChaMis} consists in first proving it for countable acyclic groups, and
then embedding an arbitrary countable group $G$ in a countable acyclic group,
such as $A(G)$ (using (d) of Theorem \ref{Berrick} above).

\appendix

\section{Examples of low-dimensional acyclic groups}

Examples described in \cite{Berrick_Steer chapter} that are of low
cohomological dimension are as follows.

First, Higman's four-generator, four-relator group ($k=4$)
\[
\left\langle x_{i}\mid x_{i+1}=[x_{i},x_{i+1}]\right\rangle _{i\in
\mathbb{Z}/k}%
\]
(with $k\geq0$, but $k=1,2,3$ trivial) is a candidate for being the
\textquotedblleft oldest\textquotedblright\ acyclic group in the literature,
although its acyclicity was not proved until much later in \cite{DyerVasquez},
where it is shown to be the fundamental group of both an acyclic, aspherical
finite $2$-complex and an homology $4$-sphere. In \cite{BerrickWong} it is
observed that this group is the commutator subgroup of
\[
\left\langle x,y\mid\,x[x,yxy^{-1}],\ [x,y^{4}]\right\rangle \text{.}%
\]
For a discussion of how the Higman group relates torsionfree acyclic groups
and groups $G$ with \underline{$E$}$G/G$ contractible, see \cite{LearyNucinkis
Every complex classifying space for proper bundlesTopology 01}.

Stitch groups \cite{BerrickWong} are generalizations of the wild arc groups of
\cite{FoxArtin}, corresponding to iteration (infinitely, towards a cluster
point, in both directions) of an underlying $2$-tangle called a \emph{stitch}.
They are also acyclic of Baumslag-Gruenberg type. In fact, with the
\emph{deficiency} $\mathrm{def}(G)$ as the maximum excess of the number of
generators over the number of relators in a finite presentation of $G$, the
following hold.

\begin{theorem}
\label{01228a}\label{00d19a}\cite{BerrickHillman}, \cite{BerrickWong} Let
$\lambda$ be a smooth knot in $S^{2}\times S^{1}$ such that $[\lambda]$
generates $H_{1}(S^{2}\times S^{1})$, and let $M$ be the closed complement of
a tubular neighbourhood of $\lambda$ in $S^{2}\times S^{1}$. Then

\begin{enumerate}
\item[(a)] $M$ is an aspherical Haken manifold, and is an integral homology circle;

\item[(b)] $\pi=\pi_{1}(M)$ is finitely presentable, $\mathrm{def}(\pi)=1$ and
$\mathrm{gd}(\pi)\leq2$;

\item[(c)] $\pi$ is locally indicable and in particular contains no nontrivial
finitely generated perfect subgroups;

\item[(d)] $\pi$ is residually finite;

\item[(e)] the commutator subgroup $\pi^{\prime}$ is acyclic.
\end{enumerate}
\end{theorem}

Here is an algebraic counterpart of this situation.

\begin{theorem}
\label{00d30b}\cite{BerrickHillman} Let $G$ be a finitely presentable group
with $\mathrm{def}(G)>0$, and such that $G^{\prime}$ is finitely generated and
perfect, nontrivial. Then the following hold.

\begin{enumerate}
\item[(a)] $\mathrm{def}(G)=1$, $G$ requires at least three generators, and
$\mathrm{gd}(G)\leq2$.

\item[(b)] $G_{\mathrm{ab}}$ is $\mathbb{Z}$ or $\mathbb{Z}^{2}$.

\item[(c)] $G^{\prime}$ is acyclic, but not $FP_{2}$.
\end{enumerate}
\end{theorem}

\section{Examples of binate groups}

Examples discussed in \cite{Berrick_Steer chapter} include:

\begin{itemize}
\item algebraically closed groups (which are relevant to the study of the
Kervaire conjecture cited above \cite{BerrVarad});

\item Philip Hall's countable universal locally finite group;

\item groups of self-homeomorphisms with support symmetries;

\item the general linear group of the cone on a ring;

\item the cone on a group $G$, as constructed in \cite{KanThurston} as the
semidirect product of $G^{\mathbb{Q}}$ by $\mathrm{Aut}(\mathbb{Q})$. Here,
$G^{\mathbb{Q}}$ is the set of functions from the rationals $\mathbb{Q}$ to
$G$ which map all numbers outside some finite interval to the identity; the
group structure on $G$ determines that on $G^{\mathbb{Q}}$. Likewise,
$\mathrm{Aut}(\mathbb{Q})$ denotes the restricted symmetric group on
$\mathbb{Q}$ comprising those permutations with compact support. Here, $G$ is
embedded as a two-step subnormal subgroup of a binate group.

\item mitotic groups of \cite{BDH} providing combinatorial, finitely generated
and, later \cite{BDM}, finitely presented results analogous to those of
\cite{KanThurston}. They have the remarkable property that all quotients are
also binate.

\item automorphism groups of large structures elaborated in \cite{delaHMcD}.
Examples include the group of all continuous linear automorphisms, or of
invertible isometries, of an infinite-dimensional Hilbert space, the group of
invertible or of unitary elements in a properly infinite von Neumann algebra,
the group of measure-preserving automorphisms of a Lebesgue measure space, and
the group of permutations on an arbitrary infinite set.
\end{itemize}

\bigskip

\bigskip

\noindent\textbf{Acknowledgements. }As noted above, comments of M. V. Sapir
and I. Agol have been useful to this work. The assistance of NUS Research
Grant R-146-000-097-112 is gratefully acknowledged.

\end{document}